\documentclass[a4paper,12pt]{amsart}
\usepackage{amssymb}
\usepackage{booktabs}
\usepackage{color}
\usepackage{hyperref}
\usepackage{nicefrac}

\usepackage{tikz}
\usetikzlibrary{patterns}
\usetikzlibrary{decorations.shapes}
\usepackage{enumitem}

\usepackage{todonotes}

\newtheorem{theorem}{Theorem}[section]
\newtheorem{lemma}[theorem]{Lemma}
\newtheorem{proposition}[theorem]{Proposition}

\newtheorem{corollary}[theorem]{Corollary}
{\theoremstyle{remark}
\newtheorem{remark}[theorem]{Remark}
\newtheorem{question}[theorem]{Question}

}
\theoremstyle{definition}

\newtheorem{example}[theorem]{Example}

\newcommand{\Amp}{\mathrm{Amp}}
\newcommand{\Stab}{\mathrm{Stab}}
\newcommand{\sStab}{\mathrm{sStab}}
\newcommand{\E}{\mathcal{E}}

\newcommand{\F}{\mathcal{F}}
\newcommand{\T}{\mathcal{T}}
\renewcommand{\O}{\mathcal{O}}
\newcommand{\CO}{\mathcal{O}}

\newcommand{\FF}{\mathbb{F}}
\newcommand{\RR}{\mathbb{R}}
\newcommand{\QQ}{\mathbb{Q}}

\newcommand{\ZZ}{\mathbb{Z}}
\newcommand{\NN}{\mathbb{N}}
\newcommand{\PP}{\mathbb{P}}

\newcommand{\KK}{k}

\DeclareMathOperator{\vol}{vol}

\DeclareMathOperator{\rk}{rk}
\DeclareMathOperator{\Span}{span}

\DeclareMathOperator{\Nef}{Nef}
\DeclareMathOperator{\Mov}{Mov}
\newcommand{\Hom}{\mathcal{H}om}

\DeclareMathOperator{\rank}{rank}

\DeclareMathOperator{\conv}{conv}

\newcommand{\lin}{\mathrm{lin}}
\newcommand{\cdits}{\cdots\stackrel{\stackrel{i}{\vee}}{\vphantom{a}}\cdots}
\newcommand{\rotateRPY}[4][0/0/0]% point to be saved to \savedxyz, roll, pitch, yaw
{   \pgfmathsetmacro{\rollangle}{#2}
    \pgfmathsetmacro{\pitchangle}{#3}
    \pgfmathsetmacro{\yawangle}{#4}

    \pgfmathsetmacro{\newxx}{cos(\yawangle)*cos(\pitchangle)}% a
    \pgfmathsetmacro{\newxy}{sin(\yawangle)*cos(\pitchangle)}% d
    \pgfmathsetmacro{\newxz}{-sin(\pitchangle)}% g
    \path (\newxx,\newxy,\newxz);
    \pgfgetlastxy{\nxx}{\nxy};

    \pgfmathsetmacro{\newyx}{cos(\yawangle)*sin(\pitchangle)*sin(\rollangle)-sin(\yawangle)*cos(\rollangle)}% b
    \pgfmathsetmacro{\newyy}{sin(\yawangle)*sin(\pitchangle)*sin(\rollangle)+ cos(\yawangle)*cos(\rollangle)}% e
    \pgfmathsetmacro{\newyz}{cos(\pitchangle)*sin(\rollangle)}% h
    \path (\newyx,\newyy,\newyz);
    \pgfgetlastxy{\nyx}{\nyy};

    \pgfmathsetmacro{\newzx}{cos(\yawangle)*sin(\pitchangle)*cos(\rollangle)+ sin(\yawangle)*sin(\rollangle)}
    \pgfmathsetmacro{\newzy}{sin(\yawangle)*sin(\pitchangle)*cos(\rollangle)-cos(\yawangle)*sin(\rollangle)}
    \pgfmathsetmacro{\newzz}{cos(\pitchangle)*cos(\rollangle)}
    \path (\newzx,\newzy,\newzz);
    \pgfgetlastxy{\nzx}{\nzy};

    \foreach \x/\y/\z in {#1}
    {   \pgfmathsetmacro{\transformedx}{\x*\newxx+\y*\newyx+\z*\newzx}
        \pgfmathsetmacro{\transformedy}{\x*\newxy+\y*\newyy+\z*\newzy}
        \pgfmathsetmacro{\transformedz}{\x*\newxz+\y*\newyz+\z*\newzz}

    }
}

\tikzset{RPY/.style={x={(\nxx,\nxy)},y={(\nyx,\nyy)},z={(\nzx,\nzy)}}}

\title[Stability of tangent bundles on toric surfaces...]{Stability of tangent bundles on smooth toric Picard-rank-$2$ varieties and surfaces}
\author[M. Hering]{Milena Hering}
\address{Milena Hering, School of Mathematics, University of Edinburgh, Peter-Guthrie-Tait Road, EH9 3EF, Edinburgh, United Kingdom}
\email{\href{mailto:m.hering@ed.ac.uk}{m.hering@ed.ac.uk}}
\author[B. Nill]{Benjamin Nill}
\address{
Benjamin Nill, Fakult\"at f\"ur Mathematik, 
Otto-von-Guericke-Universit\"at Magdeburg, 
Universit\"atsplatz 2, 
39106 Magdeburg, Germany
}
\email{\href{mailto:benjamin.nill@ovgu.de}{benjamin.nill@ovgu.de}}
\author[H. S\"u{\ss}]{Hendrik S\"u\ss}
\address{
Department of Mathematics,
The University of Manchester,
Alan Turing Building,
Oxford Road,
Manchester M13 9PL, United Kingdom 
}
\email{\href{mailto:hendrik.suess@manchester.ac.uk}{hendrik.suess@manchester.ac.uk}}

\begin{document}
\begin{abstract}
 We give a combinatorial criterion for the tangent bundle on a smooth toric variety to be stable with respect to a given polarisation in terms of the corresponding lattice polytope. Furthermore, we show that for a  smooth toric surface $X$ and a smooth toric variety of Picard rank 2,  there exists an ample line bundle with respect to which the tangent bundle is stable if and only if it is an iterated blow-up of projective space. 
\end{abstract}
\maketitle

\section{Introduction}

Let $X$ be a smooth toric variety over a field of characteristic 0, with  tangent bundle $\T_X$. Let $\O(D)$ be an ample line bundle. Recall that the slope of a torsion-free sheaf $\E$ on a normal projective variety $X$ with respect to a nef line bundle $\O(D)$ is defined to be 
\[\mu(\E)= \frac{c_1(\E)\cdot D^{n-1}}{\rank(\E)},\] and that $\E$ is 
is \emph{stable} (resp. \emph{semistable}) with respect to $\O(D)$ if for any subsheaf $\F$ of $\E$ of smaller rank, we have $\mu(\F) < \mu(\E)$ (resp. $\mu(\F)\leq \mu(\E)$). A direct sum of stable sheaves with identical slope is called \emph{polystable}. A situation of particular interest is when $X$ is Fano, $\E = \T_X$ is the tangent bundle, and $D = -K_X$ the anticanonical divisor, in particular, since the existence of a K{\"a}hler-Einstein metric on a Fano variety implies that the tangent bundle is polystable with respect to the anticanonical polarisation, see Section~\ref{sec:KE} for more details.

The main question we discuss in this article is when toric varieties admit a polarisation $\O(D)$ such that the tangent bundle $\T_X$ is stable with respect to $\O(D)$. This question has been studied in \cite{Pang} and recently also by Biswas, Dey, Genc, and Poddar in \cite{BiswasDeyGenvPoddar}. Note that it is well-known that the tangent bundle on projective space is stable with respect to $\mathcal{O}_{\PP^n}(1)$.

\begin{theorem}\label{thm:Easy}
Let $X$ be a smooth toric surface or a smooth toric variety of Picard rank~2. Then there exists an ample line bundle $\mathcal{L}$ on $X$ such that $T_X$ is stable with respect to $\mathcal{L}$ if and only if it is an iterated blow-up of  projective space. 
\end{theorem}

For more precise statements, see Theorems~\ref{thm:surface-stability} and \ref{thm:kleinschmidt}.
Theorem~\ref{thm:kleinschmidt} and a more detailed discussion of the Fano case has been independently obtained by Dasgupta, Dey, and Khan \cite{DDK}.
While for smooth toric varieties of Picard rank 3 it is an open question whether Theorem~\ref{thm:Easy} holds, there exists a toric Fano 3-fold of Picard rank 4 whose tangent bundle is stable with respect to the anticanonical polarisation, but that does not admit a morphism to $\PP^3$, see Example~\ref{ex:NoMapsToPn}.

 We deduce the following criterion for the tangent bundle $\T_X$ on a toric variety $X$ to be stable with respect to a given polarisation $\O(D)$ from well-known descriptions of stability conditions 
 in terms of the Klyachko filtrations associated to the tangent bundle (see, for example, \cite{KlyachkoHermitian,zbMATH01324687,zbMATH05931673}).
Let $P_D$ be the lattice  polytope associated to $D$. For each ray $\rho$ in the fan $\Sigma$, let $P_D^{\rho}$ denote the facet corresponding to $\rho$. 
\begin{proposition}
\label{thm:tangent-criterion}
	The tangent bundle on a smooth projective toric variety $X$ of dimension~$n$ is (semi)-stable with respect to an ample line bundle $\O(D)$ on X if and only if for every proper subspace $F \subsetneq N \otimes k$ the following inequality holds:
\begin{equation}
\label{eq:tangent-stability}
    \frac{1}{\dim F} \sum_{v_\rho  \in F} \vol(P_D^\rho) \stackrel{(\leq)}{<}  \frac{1}{n} \sum_{\rho} \vol(P_D^\rho) =: \frac{1}{n} \vol \partial P_D.
\end{equation}

\end{proposition}
Here, $\vol(P^\rho)$ denotes the lattice volume inside the affine span of $P^\rho$ with respect to the lattice $\Span(P^\rho) \cap M$.

We now present our results with more details.  
Let $\Amp(X) \subset N^1(X)_{\mathbb{R}}$ denote the ample cone of $X$.
It is convenient to define 
\begin{eqnarray*} \Stab(\T_X) &=& \{ D \in \Amp(X) \mid \T_X \text{ is stable with respect to } \CO(D)\} \textrm{, and} \\ 
\sStab(\T_X)& =& \{ D \in \Amp(X) \mid \T_X \text{ is semistable with respect to } \CO(D)\}.\end{eqnarray*}

Using results from \cite{zbMATH06561543} one can show that if for a $\QQ$-factorial variety $\Stab(\T_X)$ is non-empty, then for any birational morphism $X'\to X$, $\Stab(\T_{X'})\neq \emptyset$, see \ref{lem:blow-up}. In particular, since the tangent bundle to $\PP^n$ is stable with respect to the anticanonical polarisation, 
any iterated blow-up of projective space has $\Stab(\T_X) \neq \emptyset$. 

Recall that every smooth toric surface is either a successive toric blow-up of $\PP^2$ or of  a Hirzebruch surface $\FF_a$. In Lemma~\ref{lem:not-blowup-of-p2}, we characterise the fans of smooth toric surfaces that are not a blow-up of $\PP^2$ and use this to prove the following theorem. 

\begin{theorem}
  \label{thm:surface-stability}
Let $X$ be a smooth toric surface. Then  
	\begin{enumerate}
		\item $\Stab(\T_X) = \Amp(X)$ if and only if  $X = \PP^2$
                \item $\emptyset = \Stab(\T_X) \subsetneq \sStab(\T_X)$ if and only if $X \cong \PP^1 \times \PP^1$.
		\item $\emptyset \subsetneq \Stab(\T_X) \subsetneq \Amp(X)$ if and only if $X$ is an iterated blow-up of $\PP^2$, but not $\PP^2$ itself, \label{item:proper-blow-up}
		\item $\Stab(\T_X) = \emptyset$  if and only if $X$ is not an iterated blow-up of $\PP^2$. 
	\end{enumerate}
\end{theorem}

In  \cite[Theorem 6.2]{BiswasDeyGenvPoddar}, Biswas et al. show that when $X$ is the Hirzebruch surface $\FF_a$, $a\geq 2$ implies that $\Stab(\T_X)=\emptyset$ and for $a=1$ they describe $\Stab(\T_X)$
in \cite[Corollary 6.3]{BiswasDeyGenvPoddar}. 

Projectivisations of direct sums of line bundles on projective spaces yield examples of toric Fano varieties under some conditions, but are also interesting in their own right. 
By \cite[Theorem 1]{zbMATH04085806} every smooth toric variety of Picard rank 2 is of the form $X = \PP_{\PP^s}(\CO \oplus \bigoplus_{i=1}^r \CO (a_i))$, and $X$ is a blow-up of $\PP^s$ if and only if $(a_1, \ldots, a_r) = (0, \ldots, 0, 1)$.
Note that the polytopes corresponding to ample line bundles on these varieties  are special cases of \emph{Cayley polytopes}, see for example \cite{CattaniCoxDickenstein97}.

\begin{theorem}
\label{thm:kleinschmidt}
Consider the smooth projective variety
\[X = \PP_{\PP^s}(\CO \oplus \bigoplus_{i=1}^r \CO (a_i))\]
for $s,r \ge 1$ with $0 \le a_1 \le \cdots \le a_r$. For $a_r \ge 1$, we have $\Stab(\T_X) \neq \emptyset$ if and only if $\sStab(\T_X) \neq \emptyset$ if and only if $(a_1, \ldots, a_r) = (0, \ldots, 0, 1)$. In this case, $\T_X$ is (semi-)stable with respect to a polarisation $\mathcal{L} = \CO_X(\lambda) \otimes \pi^*\CO(\mu)$ if and only if $p(\mu/\lambda) \stackrel{(\le)}{<} 0$, where $p(x)$ is the following polynomial of degree $s$:
\[p(x) := -\left(\sum_{q = 0}^{s-1} \binom{r+s-1}{q} x^{q}\right) + \frac{s(r+1)}{r} \binom{r+s-1}{s} x^{s}.\]
	We note that $p(\mu/\lambda) < 0$ if and only if $\mu/\lambda$ is in the interval $(0,\gamma)$, where $\gamma$ is the only positive root of $p(x)$. For $r=1$ we have $\gamma = \frac{1}{(2s+1)^{\nicefrac{1}{s}}-1}$, and for $s=1$ we get $\gamma= \frac{1}{r+1}$.

One has $\emptyset = \Stab(\T_X) \subsetneq \sStab(\T_X)$ if and only if $(a_1,\ldots, a_r)=(0,\ldots,0)$, i.e. if $X = \PP^s \times \PP^r$. In this case $\T_X$ is semistable only with respect to pluri-anticanonical polarisations.
\end{theorem}

This result has been independently proved by \cite{DDK}. 
It is extending a result by Biswas et. al. \cite[Theorem 8.1]{BiswasDeyGenvPoddar}, who show that in the  Fano case (when $0<a\leq s$),  and when $s\geq 2$, the tangent bundle on $X=\PP_{\PP^s}(\CO \oplus \CO(a))$ is not stable with respect to the anticanonical polarisation $\CO(-K_X) = \CO(2)\otimes \pi^{*}\CO(s+1-a)$.

The tangent bundle to a smooth Fano surface is stable with respect to the anticanonical polarisation by \cite{zbMATH04099463}. Moreoever, in \cite{zbMATH00892566} all smooth  Fano threefolds with stable (resp. semistable) tangent bundle are classified. 
 Moreover, for smooth toric Fano varieties of dimension 4 and Picard rank 2,  the (semi-)stability of the tangent bundle with respect to the anticanonical polarisation is treated in \cite[Section 9]{BiswasDeyGenvPoddar}, and  for smooth toric Fano varieties of dimension 4 and Picard rank 3 in \cite{DDK}.

The above results motivate the following question:

\begin{question}Are there only finitely many isomorphism classes of smooth projective toric varieties $X$ of given dimension $n$ and Picard number $\rho$ with $\Stab(T_X) \not=\emptyset$?
\label{question-finite}
\end{question}

\begin{corollary}Question~\ref{question-finite} has an affirmative answer for $n\le 2$ or $\rho\le 2$.
\end{corollary}

\begin{proof}
The cases $n=1$ or $\rho=1$ are trivial. For $n=2$ this follows from Theorem~\ref{thm:surface-stability}(3). For $\rho=2$ this follows from Theorem~\ref{thm:kleinschmidt} (note that $\dim(X) = r+s$).
\end{proof}

\subsection{Connections to the existence problem of K\"ahler-Einstein metrics}\label{sec:KE}
When $X$  is a smooth Fano variety over the complex numbers, the existence of a K{\"a}hler-Einstein metric on the underlying complex manifold $X$ implies that its tangent bundle is polystable, (in particular, semistable) with respect to the anticanonical polarisation \cite{zbMATH03890301}, \cite[Sec 5.8]{zbMATH00044936}. However, the converse does not hold for the blow-up of $\PP^2$ in two points \cite{zbMATH03149530}. 
The recent proof of the Yau-Tian-Donaldson conjecture \cite{zbMATH06497829,zbMATH06394344,zbMATH06394345,zbMATH06394346} shows that a Fano manifold has a K{\"a}hler-Einstein metric if and only if it is $K$-polystable. For a general toric Fano variety K-stability is equivalent to the fact that for the polytope corresponding to the anticanonical polarisation the barycenter coincides with the origin \cite{li2018algebraicity}, in the smooth case this was known before due to combining  \cite{wang04} and \cite{Mabuchi87}.  

Thus we obtain the following combinatorial statement: 
\begin{corollary}
	Let $P$ be a smooth reflexive polytope with barycenter in the origin. Then $P$ satisfies the non-strict inequality
 (\ref{eq:tangent-stability}) for every proper linear subspace $F \subset N_\QQ$. 
\end{corollary}

This statement has been known to combinatorialists in a more general setting that implies the statement for reflexive polytopes with barycenter in the origin (without the smoothness assumption). 
Conditions of this type are known in convex geometry under the name \emph{subspace concentration conditions}. They play a distinguished role in several problems from convex geometry, see e.g. \cite{zbMATH06435124,zbMATH06168607,zbMATH06285042}. The fact that this condition holds for a reflexive polytope whenever the barycenter coincides with the origin is far from being obvious. Moreover, our argument via K\"ahler-Einstein metrics is valid only in the smooth case (since we have to rely on \cite{zbMATH03890301}, \cite[Sec 5.8]{zbMATH00044936}), but the fact turns out to be true for every reflexive polytope. This follows from an even more general result in \cite[Thm~1.1]{zbMATH06285042}, which applies to every polytope with barycentre at the origin. Their proof relies entirely on methods from convex geometry.

\section*{Acknowledgements} 
We would like to thank 
Carolina Araujo, Arend Bayer, Stefan Kebekus, and Adrian Langer for helpful discusstions.  We would also like to thank
University of Edinburgh, Institute of Mathematics of the Polish Academy of Sciences  Warsaw,  and Mathematical Sciences Research Institute for their hospitality. 
BN is funded by the Deutsche Forschungsgemeinschaft (DFG, German Research Foundation) - 314838170, GRK 2297 MathCoRe. BN is also partially supported by the Vetenskapsr\aa det grant NT:2014-3991 (as an affiliated researcher with Stockholm University).

\section{Stability conditions for equivariant sheaves}
\label{sec:stab-cond-toric}
We fix our setting as follows. We consider a polarized toric variety $(X, \CO(D))$ corresponding to a lattice polytope $P$. Let $\Sigma$ be the normal fan of $P$ and $P^\rho$ the facet of $P$ corresponding to a ray $\rho \in \Sigma$.

Recall that a coherent sheaf $\E$ is called {\em reflexive} if $\E \cong \E^{\vee\vee}$,
where $\E^\vee = \Hom(\E,\CO_X)$. In \cite{Klyachko90} equivariant vector bundles on smooth toric varieties were classified in terms of collections of filtrations of $k$-vector spaces indexed by the rays of $\Sigma$. This classficiation extends to equivariant reflexive sheaves on normal toric varieties, see for example, \cite{zbMATH01324687, Perling04}. 

More precisely, we fix a $k$-vector space $E$ and for every ray $\rho \in \Sigma^{(1)}$ we consider a decreasing filtration by subspaces
\[E \supset \ldots \supset E^{\rho}(i-1) \supset E^{\rho}(i) \supset E^{\rho}(i+1) \supset \ldots \supset 0,\]
such that $E^{\rho}(i)$ differs from $E$ and $0$ only for finitely many values of $i \in \ZZ$.
Given such a collection of filtrations for every cone $\sigma \subset \Sigma$ we may consider
\[
 E_u:=\left(\bigcap_{\rho\in \Sigma^{(1)}} E^\rho(-\langle v_\rho, u \rangle)\right) \otimes \chi^u \subset E \otimes k[M].
\]
Then $\bigoplus_{u \in M} E_u$ is equipped with the structure of an $M$-graded $k[U_\sigma]$-module via the natural multiplication with $\chi^u \in k[U_\sigma]$. Then setting $H^0(U_\sigma,\E)=\bigoplus_{u \in M} E_u$
for every $\sigma \in \Sigma$ defines an equivariant reflexive sheaf on $X$.

The collections of filtrations form an abelian category in a natural way. A morphism between a collection of filtrations $F^{\rho}(i)$ of a vector space $F$ and another collection $E^{\rho}(i)$ of filtrations of a vector space $E$ a morphism is a linear map $L \colon F \to E$ which are compatible with the filtrations, i.e. 
$L(F^\rho(i)) \subset E^\rho(i)$ for all $\rho \in \Sigma^{(1)}$ and all $i \in \ZZ$.

\begin{theorem}
  There is an equivalence of categories between the equivariant reflexive sheaves on a toric variety $X = X_\Sigma$ and the collections of filtrations of $k$-vector spaces indexed by the rays of $\Sigma$. Here, the rank of the reflexive sheaf equals the dimension of the filtered $k$-vector space.
\end{theorem}

For a collection of filtrations $E^{\rho}(i)$,  we set $e^{[\rho]}(i)=\dim E^{\rho}(i) - \dim E^{\rho}(i+1)$ similarly for other filtrations we will always use the lower letter version to denote the differences of dimensions between the steps of the filtration. Then we have the following formula.

\begin{lemma}
  \label{lem:mu}
  Assume that $X$ is smooth. With the notation above we have
  \[
    \mu(\E) = \frac{1}{\dim E} \sum_{i,\rho} i \cdot e^{[\rho]}(i) \cdot \vol(P^\rho).
  \]
\end{lemma}

\begin{proof}
  By \cite[Corollary 3.18]{zbMATH05931673},  we have
  $c_1(\E)=\sum_{\rho} \sum_{i \in \ZZ} i  e^{[\rho]}(i)  D_\rho$. Now,
  for a ray $\rho \in \Sigma^{(1)}$ the intersection number
  $D^{n-1} \cdot D_\rho$ is given by the the volume of the corresponding facet $P^\rho$ of $P$, see e.g. \cite{zbMATH03661497}.
\end{proof}

With the notation above we get the following characterisation of stability.
\begin{proposition}
\label{prop:saturated-subsheaf-criterion}
Let $X$ be a smooth toric variety. A toric vector bundle $\E$ on $X$ corresponding to filtrations $E^\rho(i)$ is (semi-)stable if and only if the following  inequality holds for every linear subspace $F \subset E$ and $F^\rho(i)=E^\rho(i) \cap F$.
        \begin{equation}
          \label{eq:stability}
          \frac{1}{\dim F} \sum_{i,\rho} i \cdot f^{[\rho]}(i) \cdot \vol(P^\rho) \stackrel{(\leq)}{<} \frac{1}{\dim E} \sum_{i,\rho} i \cdot e^{[\rho]}(i) \cdot \vol(P^\rho)
        \end{equation}
\end{proposition}
\begin{proof}

  By \cite[Proposition 4.13]{zbMATH05931673} it is sufficient to consider equivariant reflexive subsheaves. It remains to show that it is sufficient to consider those subsheaves, which correspond to filtrations of the form $E^{\rho}(i) \cap F$. For every subsheaf $\F' \subset \E$ corresponding to filtrations $(F')^\rho(i) \subset E^\rho(i)$ of some subspace $F \subset E$ we may consider the subsheaf $\F$ corresponding to the filtrations $F^{\rho}(i):=E^{\rho}(i) \cap F$. Then $\dim F^\rho(i) \geq \dim (F')^\rho(i)$ for all $i,\rho$. Now,  Lemma~\ref{lem:ineq} implies that $\mu(\F) \geq \mu(F')$.
\end{proof}

\begin{remark}
  A subsheaf $\F$ of a torsion-free sheaf $\E$ is called {\em saturated} if $\E/\F$ is torsion-free.  The saturation of a subsheaf $\F \subset \E$ is the smallest saturated subsheaf of $\E$ containing $\F$. It is not hard to derive from the description of $H(U_\sigma,\E)$ given above, that $\F \subset \E$ given by $F^\rho(i) \subset E^\rho(i)$ is saturated, if and only if $F^\rho(i) = E^\rho(i) \cap F$. Hence, Lemma~\ref{lem:ineq} below can be seen as a combinatorial version of the well-known fact that replacing a subsheaf by its saturation increases the slope.
\end{remark}

\begin{lemma}
\label{lem:ineq}
	Given integer functions $f,g:\ZZ \to \ZZ$ with $f \geq g$ such that $\{i\in \mathbb{Z}\mid f(i)\neq g(i)\}$ is finite.  Then also 
\[\sum_i i\cdot (f(i)-f(i+1)) \geq \sum_i i\cdot (g(i)-g(i+1))\]
holds.
\end{lemma}
\begin{proof}
	Note that the assumption implies that $A(f,g):=\sum _i (f(i)-g(i)) \ge 0$ is finite. We fix $f$ and proceed by induction on $A(f,g)$. If $A(f,g)=0$, $f=g$  and the statement is trivially true. Fix $f$ and assume that the statement holds for all $g\leq f$ with $A(f,g)\leq A$. Let $g'$ be such that $A(f,g') = A+1$. Since $A> 0$, there exists a $k$ such that $f(k)>g'(k)$. Define 
	
	\[ g(i) =\begin{cases} g'(i) & \text{ if } i\neq k\\ 
		g'(i)+1 &\text{ if } i=k.
	\end{cases}\]
	Then $A(f,g) = A$. We calculate   
$\sum_i i\cdot (g(i)-g(i+1)) =  \sum_i i\cdot (g'(i)-g'(i+1)) +1$.   
	By induction hypothesis, we have  
$\sum_i i\cdot (f(i)-f(i+1)) \geq \sum_i i\cdot (g(i)-g(i+1)) >   
\sum_i i\cdot (g'(i)-g'(i+1))$.   
\end{proof}

By \cite{Klyachko90} the filtrations of the tangent bundle on $\T_X$ on a smooth toric variety $X$ have the following form.
\begin{equation}
 E^{\rho}(j) =
 \begin{cases}
   N \otimes \KK & j < 1\\
   \Span_\KK (v_\rho)& j = 1 \\
   0 & j > 1
 \end{cases}
\label{eq:TX-filtrations}
\end{equation}

\begin{proof}[Proof of Proposition~\ref{thm:tangent-criterion}]
  Looking at the filtrations $E^\rho(i)$ for $\T_X$ from (\ref{eq:TX-filtrations}) we see that
\[e^{[\rho]}(i) = E^\rho(i)-E^\rho(i+1)
=
\begin{cases}
  n-1 & j=0\\
  1   & j = 1\\
  0   & \text{ else.}
\end{cases}
\]
Similary, for  $F^\rho(i) =  E^\rho(i)\cap F$ we have
\[f^{[\rho]}(i)=
\begin{cases}
  \dim(F)-2 & j=0\\
  1   & j = 1\\
  0   & \text{else},
\end{cases}\quad
\text{ or } \quad
f^{[\rho]}(i)
=
\begin{cases}
  \dim(F)-1 & j=0\\
  0   & \text{else},
\end{cases}
\]
depending on whether $v_\rho$ is contained in the subspace $F$ or not. Now, Proposition~\ref{prop:saturated-subsheaf-criterion} immediately implies the claim of Proposition~\ref{thm:tangent-criterion}.
\end{proof}

\begin{remark}
\label{rem:finitely-many-subspaces}
  Actually, it is suffices to test the inequality of Proposition~\ref{thm:tangent-criterion} for the (finitely many) subspaces of the form $F=\Span_k R$ with $R \subset \Sigma(1)$.
Indeed, assume that $\F$, given by some $F \subset N \otimes \KK$, destabilises $\T_X$. Then we may choose $\F'$ corresponding to $F' := \Span \{v_\rho \subset \Sigma(1) \mid v_\rho \subset F\}$. With this choice
we have $\sum_{v_\rho \in  F} \vol(P^\rho) = \sum_{v_\rho \in F'} \vol(P^\rho)$ and
$\rk \F' =\dim F' \leq \dim F = \rk \F$.
\end{remark}

\begin{example}
\label{exp:ppn}
  For $\PP^n$ a polarisation is given by $\CO(d)$. The corresponding polytope is a $d$-fold dilation of the standard simplex $\Delta \subset \RR^n$. Every facet of $d\Delta$ has lattice volume $d^{n-1}$. For every proper subset $R \subsetneq \Sigma(1)$ and $F= \Span R$  we have $\dim F = \# R$. Now (\ref{eq:tangent-stability}) becomes $d^{n-1} < d^{n-1} \cdot (n+1)/n$. Hence, we recover the well-known fact, that $\PP^n$ has a stable tangent bundle.
\end{example}

\begin{lemma}
  \label{lem:blow-up}
Assume that $X$ is $\QQ$-factorial and $\Stab(\T_X)$ is non-empty. If there is a birational morphism $f \colon X' \to X$ , then $\Stab(\T_{X'})$ is non-empty, as well.
\end{lemma}
\begin{proof}
Consider a polarisation $\CO(D)$ of $X$, such that $\T_X$ is stable. Then $\T_{X'}$ is stable with respect to the nef and big bundle $\CO(f^*D)$, since any destabilising subsheaf $\F' \subset \T_{X'}$ with respect to $\CO(f^*D)$ would induce a subsheaf $(f_*\F') \subset \T_X$  which, by projection formula, would be destabilising with respect to $\CO(D)$. Now, the openness property from \cite[Thm 3.3]{zbMATH06561543} ensures the existence of a stabilising ample class, which is given as a small pertubation of $[\CO(f^*D)]$.
\end{proof}

We also have the following equivalent for the strictly unstable case.
\begin{lemma}
  \label{lem:blow-up-unstable}
Assume that $X$ is $\QQ$-factorial and $\Amp(X) \setminus \sStab(\T_X)$ is non-empty. If there is a birational morphism $f \colon X' \to X$ , then $\Amp(X') \setminus \sStab(\T_{X'})$ is also non-empty.
\end{lemma}
\begin{proof}
 Assume that a subsheaf $\F \subset \T_X$ destabilises $\T_X$ strictly with respect to an ample polarisation $\mathcal{O}(D)$. Then we 
 note that $f^*\F$ and $\T_{X'}$ are both subsheaves of $f^*\T_X$. Now, we claim that $\F':=f^*\F \cap \T_{X'}$ destabilises $\T_{X'}$ with respect 
 to $\mathcal{O}(D')=f^*\mathcal{O}(D)$. Indeed, by the projection formula we obtain
 \begin{align*}
 c_1(\T_{X'}).(D')^{n-1} &=  c_1(\T_{X}).(D)^{n-1} \\
 c_1(\F').(D')^{n-1} &=  c_1(\F).(D)^{n-1}.
 \end{align*}
 The line bundle $\mathcal{O}(D')$ is only nef, but the condition that a subsheaf destabilises strictly
 is an open condition on the divisor class.  Hence, we can find an ample divisor class with the same property as a small perturbation of $\mathcal{O}(D')$.
 \end{proof}

\section{Smooth toric surfaces}
Every toric surface can be obtained via equivariant blow-ups from $\PP^2$ or from a Hirzebruch surface $F_a=\PP_{\PP^1}(\CO_{\PP^1} \oplus \CO_{\PP^1}(a))$, see e.g. \cite{zbMATH00192947}. For $\PP^2$ it is well-known that the tangent bundle is stable (see also Example~\ref{exp:ppn}). The following corollary, which can be also found e.g. in \cite[Sec.~6]{BiswasDeyGenvPoddar}, clarifies the situation for the Hirzebruch surfaces.
\begin{corollary}
\label{cor:hirzebruch}
  For a Hirzebruch surface $F_a=\PP_{\PP^1}(\CO_{\PP^1} \oplus \CO_{\PP^1}(a))$  the tangent bundle is semistable with respect to  $\CO_{F_a}(\lambda) \otimes \pi^*\CO_{\PP^1}(\mu)$ in the following cases
  \begin{enumerate}
  \item $a=0$ and $\lambda=\mu$,
  \item $a=1$ and $2\mu \leq \lambda$.
  \end{enumerate}
  The tangent bundle is stable if and only if $a=1$ and $2\mu < \lambda$.
\end{corollary}
\begin{proof}
  The claim follows directly from Proposition~\ref{thm:kleinschmidt} for the case $r=s=1$.
\end{proof}

\begin{lemma}
\label{lem:not-blowup-of-p2}
  A a smooth toric surfaces $X=X_\Sigma$ is not a blowup of $\PP^2$ or $\PP^1 \times \PP^1$ if and only if there are integers $a,c,e$ fulfilling $a \geq c > e+1 \geq 1$ such that after appropriate choice of basis for $N$ 
  \begin{enumerate}
  \item $\Sigma$ contains the rays spanned by \label{item:rays}
    \[(0,1), (1,0), \quad (0,-1), (1,-e), \quad (-1,c), (-1,a)\]
  \item all other rays are contained in the cones
$\langle (-1,c), (-1,a)\rangle$ and $\langle (1,0), (1,-e) \rangle$.
  \end{enumerate}
\end{lemma}
\begin{remark}
  Note, that in Lemma~\ref{lem:not-blowup-of-p2} we explicitly allow the cases $(1,-e)=(1,0)$
   $(-1,c) = (-1,a)$.
\end{remark}

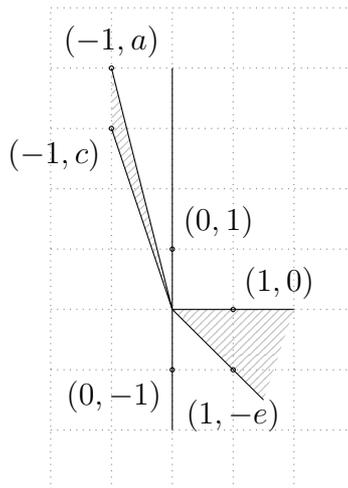
\begin{figure}[ht]
  \centering
  \begin{tikzpicture}[scale=0.8]
    \draw[dotted,step=1,gray] (-2,-3) grid (3,5);
    \draw (0,4) -- (0,-2);
    \path[pattern=north east lines,draw=black,pattern color=lightgray] (-1,4)--(0,0)--(-1,3);
    \path[pattern=north east lines, draw=black, pattern color=lightgray] (1.5,-1.5)--(0,0)--(2,0);
    \draw  (0,1) node[anchor=south west] {$(0,1)$} circle (1pt);
    \draw  (0,-1) node[anchor=north east] {$(0,-1)$} circle (1pt);
    \draw  (1,0) node[anchor=south west] {$(1,0)$} circle (1pt);
    \draw  (1,-1) node[below=7pt] {$(1,-e)$} circle (1pt);
    \draw  (-1,4) node[anchor=south] {$(-1,a)$} circle (1pt);
    \draw  (-1,3) node[anchor=north east] {$(-1,c)$} circle (1pt);
  \end{tikzpicture}

\caption{Schematic picture of a fan of a toric surface blowing down to neither $\PP^2$ nor $\PP^1 \times \PP^1$. All additional rays have to be contained in the shaded regions.}
\label{fig:fan}
\end{figure}

\begin{proof}
  For example by \cite[Thm. 1.28]{zbMATH00192947} we may find a ray $\langle v \rangle \in \Sigma(1)$ such that $-v$ spans another ray in $\Sigma(1)$. We then may number the ray generators $v_0, \ldots, v_r$ of the rays in  $\Sigma(1)$ in consecutive order, such that $v_1 = v$ and $v_\ell=-v$ for some $\ell \in \{3,\ldots, r\}$. Then by our smoothness condition $v_1,v_2$ form a basis of $N$. Hence, we may assume that $v_1=(0,1)$ and $v_2=(1,0)$. 
Again by smoothness $v_0 = (-1,a)$ for some $a \in \ZZ$. Since, $X$ is assumed not to be a blowup of $\PP^2$ or  $\PP^1 \times \PP^1$ we have 
\begin{equation}
  \label{eq:not-in-fan}
  \langle (-1,0)\rangle, \langle (-1,1)\rangle, \langle (-1,-1)\rangle \notin \Sigma^{(1)}.
\end{equation}
Hence, $a \neq 0,1,-1$. After possibly switching the role of $v_1$ and $v_\ell$ we may assume that $a \geq 2$. Now, assume $v_{\ell+1}= (b,c)$. The regularity of the cone $\langle (b,c), (0,-1)\rangle$ implies $b=-1$. Moreover, by \cite[Prop.~1.19]{zbMATH00192947} the smoothness of $X$ also implies that all the rays generated by vectors of the form $(-1,y)$ with $a\geq y \geq c$ have to be present in $\Sigma$. Now (\ref{eq:not-in-fan}) implies that $a \geq c \geq 2$. 

Similarly consider the ray $\langle v_{\ell-1} \rangle$ with $v_{\ell-1} = (d,-e)$. By regularity of the cone $\langle v_{\ell-1},v_\ell \rangle$ we must have $d=1$. Now, as before smoothness of $X$ implies by \cite[Prop.~1.19]{zbMATH00192947} that all the rays with generators $(1,-y)$ with $0 \leq y \leq e$ have to be contained in $\Sigma(1)$. On the other hand we must have $(1,-y) + (0,-1) \neq (1,-c)$ for all such $y$, since we assumed, that $X$ is not a blowup of $\PP^2$. Hence, $e \leq c-2$ must hold.
\end{proof}

\begin{proof}[Proof of Theorem~\ref{thm:surface-stability}]
The case of $\mathbb{P}^2$ is discussed in Example~\ref{exp:ppn}, and the 
  case of Hirzebruch surfaces including $\PP^1 \times \PP^1$ in Corollary~\ref{cor:hirzebruch}. For iterated blowups of $\PP^2$ we see by Example~\ref{exp:ppn} and Lemma~\ref{lem:blow-up} that $\Stab(\T_X)\neq \emptyset$. To show the 
  strict inclusion $\Stab(\T_X) \subsetneq \Amp(X)$ as claimed in item (\ref{item:proper-blow-up})
we refer to Corollary~\ref{cor:hirzebruch} together with Lemma~\ref{lem:blow-up-unstable}. It remains to show that in the other cases there exists a subbundle of $\T_X$ which destabilises $\T_X$. For this we may assume that $\Sigma$ has the form described in Lemma~\ref{lem:not-blowup-of-p2}. We also fix the notation of the proof of that lemma, i.e. we may choose a basis of $N$ order the primitive generators of rays in $\Sigma(1)$ clockwise, in such a way that 
\[
v_1= (0,1), v_2 =(1,0), v_{\ell-1}=(1,-e), v_\ell=(0,-1), v_{\ell+1}=(-1,c),v_0=(-1,a).
\]
In the following we show that $F=\Span v_1$ gives rise to a destabilising subbundle of $\T_X$. 

Let us denote the torus invariant prime divisors corresponding to $v_i$ by $D_i$ and the maximal cones $\langle v_i,v_{i+1} \rangle$ by $\sigma_i$. Assume $D=\sum_i a_i D_i$ is an ample divisor. Then the corresponding polytope 
\[P= \{u \in M_\RR \mid \forall_{i=0,\ldots,r}\colon \langle u , v_i \rangle \geq -a_i\}\]
 has normal fan equal to $\Sigma$. Its facets are given by
$P^{v_i} = P \cap \{u \mid \langle u, v_i \rangle=-a_i\}$ for $i=0,\ldots,r$ and its vertices by
\[u_{i} = \{u \mid \langle u, v_{i} \rangle=-a_{i}\} \cap \{u \mid \langle u, v_{i+1} \rangle=-a_{i+1}\}.\]

To prove that $F= \Span v_1$ gives rise to a destabilising subbundle by Proposition~\ref{thm:tangent-criterion} we have to show that 
\[\vol P^{v_1} + \vol P^{v_\ell} > \frac{1}{2} \sum_{i} \vol P^{v_i},\] or equivalently 
that 
\[\vol P^{v_1} + \vol P^{v_\ell} > \sum_{i\neq 1,\ell} \vol P^{v_i}.\]

Consider the trapezoid 
\[Q= \{u \in M_\RR \mid \forall_{i\in \{1,\ell-1,\ell,\ell+1\}} \colon \langle u , v_i \rangle \geq -a_i\}.\]
Then $P$ is contained in $Q$. The lattice points $u_\ell,u_{\ell+1}$ are also vertices of $Q$ and we have $Q^{v_\ell} = P^{v_\ell}$. We set $h=a_1+a_\ell$. This is the lattice distance between the two parallel facets $Q^{v_1}$ and $Q^{v_\ell}$ of the trapezoid. The two other facets of $Q$ lie on the lines
$\{\langle (1,-e), \cdot \rangle = -a_{\ell-1}\}$ and $\{\langle (-1,c), \cdot \rangle = -a_{\ell+1}\}$, respectively.

Then we have
$\vol Q^{v_\ell} = \vol Q^{v_1} + h\cdot (c-e)$ for the lengths of these facets.
Hence, 
\[\vol P^{v_\ell} = \vol Q^{v_\ell}= \vol Q^{v_1} + h\cdot (c-e) \geq \vol Q^{v_1}  + h \cdot 2 > 2h.\]
On the other hand the sum of the lattice length of the the remaining edges of $P$ is at most $2h$.
\end{proof}

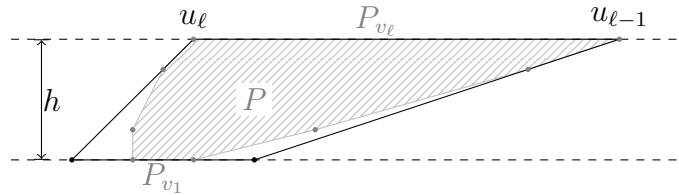
\begin{figure}[ht]
  \centering
  \begin{tikzpicture}[scale=0.8]
    \draw[draw=lightgray, pattern= north east lines,pattern color=lightgray]  (0,0)--(0,0.5)--(0.5,1.5)--(1,2)--(8,2)--(6.5,1.5)--(3,0.5)--(1,0)--cycle;
    \draw (-1,0)--(1,2)--(8,2)--(2,0)--cycle;

    \draw[dashed] (-2,0)--(9,0) ;
    \draw[dashed] (-2,2)--(9,2) ;

    \draw[fill] (-1,0) circle (1pt);
    \draw[fill] (2,0) circle (1pt);

    \begin{scope}[fill=gray,color=gray]
      \draw[fill] (0,0) circle (1pt);
      \draw[fill] (1,0) circle (1pt);
      \draw[fill] (1,2) circle (1pt);
      \draw[fill] (8,2) circle (1pt);
      \draw[fill] (0,0.5) circle (1pt);
      \draw[fill] (0.5,1.5) circle (1pt);
      \draw[fill] (3,0.5) circle (1pt);
      \draw[fill] (6.5,1.5) circle (1pt);
    \end{scope}
    \draw[|<->|] (-1.5,2) -- (-1.5,0);
    \draw (0.5,-0.3) node[color=gray] {$P_{v_1}$};
    \draw (4,2.3) node[color=gray] {$P_{v_\ell}$};
    \draw (2,1) node[color=gray,fill=white,inner sep=2pt] {$P$};
    \draw (1,2) node[anchor=south] {$u_{\ell}$};
    \draw (8,2) node[anchor=south] {$u_{\ell-1}$};
    \draw (-1,1) node[anchor=east] {$h$};
 \end{tikzpicture}

\caption{Schematic picture of the polytopes $P \subset Q$}
\label{fig:polytope}
\end{figure}

\section{Smooth toric varieties of Picard rank \texorpdfstring{$2$}{2}}
\label{sec:stable}

\subsection{The setup}
In the following we consider arbitrary projectivised vector bundles $X$ on projective spaces $\PP^s$, i.e.
\[X = \PP_{\PP^s}(\CO \oplus \bigoplus_{i=1}^r \CO (a_i))\]
for $s,r \ge 1$ and $0 \leq a_1 \leq \cdots \leq a_r$. Such a variety is Fano if and only if $\sum_i a_i \leq s$. If $a_1 = \cdots = a_r = 0$, then $X \cong \PP^s \times \PP^r$. Here, the tangent bundle splits and becomes semistable if and only if the summands have equal slope. It is straightfoward to check that this happens exactly for powers of the anti-canonical polarisation. In the following we will {\em throughout} assume that $0 \le a_1 \le \cdots \le a_r \ge 1$. Now, every ample line bundle on $X$ has the form $\mathcal{L} = \CO_X(\lambda) \otimes \pi^*\CO(\mu)$ with $\lambda, \mu > 0$.

Let us recall the notion of the Cayley sum of $(r+1)$ polytopes $P_0, \ldots, P_r$ in $\RR^s$. Consider the standard basis $e_1, \ldots, e_r$ of $\RR^r$. Then the convex hull of $(P_0 \times \{0\})$ and $(P_1 \times \{e_1\}), \ldots, (P_r \times \{e_r\})$ in $\RR^s \times \RR^r$ is denoted by $P_0 * \cdots * P_r$ and it is called the \emph{Cayley sum} of these polytopes. We remark that the polarisation by $\mathcal{L}$ chosen as above corresponds to the polytope \[\lambda (\nu \Delta_s * (a_1+\nu)\Delta_s * \cdots * (a_r+\nu) \Delta_s)) \subset M_\RR = \RR \times \RR^s,\]
where $\nu=\mu/\lambda > 0$.
Since passing to multiples of polarisations has no effect on stability, we may equivalently consider the (rational) polarisation $\CO_X(1) \otimes \pi^*\CO(\nu)$.

Hence, by Proposition~\ref{thm:tangent-criterion} we are led to investigate whether for the lattice polytope
\[P := P_D = \nu \Delta_s * (a_1+\nu)\Delta_s * \cdots * (a_r+\nu) \Delta_s\]
the subspace concentration condition \eqref{eq:tangent-stability} (strictly) holds. In this case, we say that $P$ is {\em stable} (respectively, {\em semistable}).

\subsection{The stability criterion for \texorpdfstring{$P$}{P}}

We will give in Proposition~\ref{stability-criterion} 
a criterion that allows to verify stability for $P$ without having to check condition \eqref{eq:tangent-stability} for all the subspaces spanned by primitive ray generators. For this, let us observe that the fan of $X$ (the normal fan of $P$) has the following two types of primitive ray generators in $N_\RR = \RR^r \times \RR^s$:
\begin{enumerate}[label=(\roman*)]
\item $v_0 := (-e_1 - \cdots - e_r) \times 0$, and $v_i := e_i \times 0$ for $i=1, \ldots, r$; 
\item $w_0 := (a_1 e_1 + \cdots + a_r e_r) \times (-e_1 - \cdots - e_s)$, and $w_i := 0 \times e_i$ for $i=1, \ldots, s$.
\end{enumerate}

We set $b_0 := \nu, b_1 := a_1+\nu, \ldots, b_r := a_r+\nu$. Note that $b_0 \le \cdots \le b_r$. We observe that
\begin{enumerate}[label=(\roman*)]
\item for $i=0, \ldots, r$ the facet of $P$ corresponding to $v_i$ is isomorphic to
  \[b_0 \Delta_s * \cdits * b_r \Delta_s, \]
  i.e. to the polytope obtained from the Cayley sum representation of $P$ by omitting the i-th summand.  We denote its normalized volume by $V_i$. 
\item for $i=0, \ldots, s$ the facet of $P$ corresponding to $w_i$ is isomorphic to \[b_0 \Delta_{s-1} * \cdots * b_r \Delta_{s-1}.\] We denote its normalized volume by $W$. 
\end{enumerate}

Let us note that by $b_0 \le \cdots \le b_r$, we have $V_r \le \cdots \le V_0$. Let us define the index $z \in \{0, \ldots, r-1\}$ such that 
\[0=a_1=\cdots=a_z < a_{z+1} \le \cdots \le a_r.\]

\begin{proposition} $P$ is stable (respectively, semistable) if and only if 
\[\frac{1}{r+s} (V_0 + \cdots + V_r + (s+1) W)\]
is greater than (respectively, greater than or equal to) the maximum of 
\begin{enumerate}[label=(\alph*)]
\item $V_0$
\item $\frac{1}{r} \sum_{i=0}^r V_i$
\item $W$
\item $\frac{1}{|I|+s}(\sum_{i \in I} V_i + (s+1) W)$ for $I \subseteq \{0, \ldots, r\}$ with 
$|I| < r$, and 
\begin{enumerate}[label=(\alph*)]
\item[(d1)] $\{z+1, \ldots, r\} \subseteq I$, or 
\item[(d2)] $\{0,\ldots,z\} \subseteq I$ and $|\{a_k \,:\, k \in \{z+1, \ldots, r\}\setminus I\}| = 1$.  
\end{enumerate}
\end{enumerate}
\label{stability-criterion}
\end{proposition}

For the proof we need the following observation:

\begin{lemma}
Let $I \subseteq \{0, \ldots, r\}$, and $F' :=  \lin(\{v_i \,:\, i \in I\})$.\\ Then $(a_1 e_1 + \cdots + a_r e_r) \times 0 \in F'$ if and only if 
\begin{enumerate}
\item $\{z+1, \ldots, r\} \subseteq I$, or 
\item $\{z+1, \ldots, r\} \subsetneq I$, $\{0,\ldots,z\} \subseteq I$, $|\{a_k \,:\, k \in \{z+1, \ldots, r\}\setminus I\}| = 1$.
\end{enumerate}
\label{lemma-criterion}
\end{lemma}

\begin{proof}
If (1) does not hold, then the condition $(a_1 e_1 + \cdots + a_r e_r) \times 0 \in F'$ is equivalent to the condition that one can write $a_1 e_1 + \cdots + a_r e_r = a_{z+1} e_{z+1} + \cdots + a_r e_r$ as a linear combination 
\[\lambda_0 (-e_1 - \cdots -e_r) + \sum_{i \in I \cap \{1, \ldots, z\}} \lambda_i e_i + \sum_{i \in I \cap \{z+1, \ldots, r\}} \lambda_i e_i\]
with $0 \not=\lambda_0 = \lambda_1 = \cdots = \lambda_z$, $\lambda_i = a_i + \lambda_0$ for all $i \in I \cap \{z+1, \ldots, r\}$, and  $a_j = -\lambda_0$ for all $j \in \{z+1, \ldots, r\}\setminus I$. From this the statement follows.
\end{proof}

\begin{proof}[Proof of Proposition~\ref{stability-criterion}]
Let us first show that each of the expressions appears on the left hand side of \eqref{eq:tangent-stability} . For this choose the subspace $F \subseteq \RR^r \times \RR^s$ as
(a) $F = \lin(v_0)$, (b) $F=\RR^r \times \{0\}$, (c) $F = \lin(w_0)$. In (d) note that $|I| <r$ implies that $\{v_i \;:\; i \in I\}$ is linearly independent and the spanned subspace $F' := \lin(\{v_i \;:\; i \in I\})$ does not contain any other $v_0, \ldots, v_r$. Hence, we may choose by Lemma~\ref{lemma-criterion} in (d) $F = \lin(F' \times \RR^s)$ (note that $w_0 \in F$ if and only if $(a_1 e_1 + \cdots + a_r e_r) \times 0 \in F'$). 

Now, let $F$ be any proper subspace of $\RR^r$. Clearly, we may assume that it is spanned by primitive ray generators. Let $F = \lin(\{v_i \,:\, i \in I\} \cup \{w_j \,:\, j \in J\})$ for $I \subseteq \{0, \ldots, r\}$ and $J \subseteq \{0, \ldots, s\}$. Here, we assume that $I$ and $J$ are maximally chosen, i.e., $\{v_0, \ldots, v_r\} \cap F = \{v_i \,:\, i \in I\}$ and $\{w_0, \ldots, w_r\} \cap F = \{w_j \,:\, j \in J\}$. We will distinguish two cases.

Case 1: $|I| \ge r$. In this case, we get $I = \{0, \ldots, r\}$ by maximality of $I$. Note that $\{w_0, \ldots, w_s\}$ is linearly independent. As $F$ is proper and contains $a_1 e_1 + \cdots + a_r e_r$, we get $|J| < s$ and $\dim(F)=r+|J|$. Therefore, the left hand side of \eqref{eq:tangent-stability} equals $\frac{1}{r+|J|}(V_0 + \cdots + V_r + |J| W) = \frac{1}{r+|J|}(r ((V_0 + \cdots + V_r)/r) + |J| W) \le \max((V_0 + \cdots + V_r)/r, W)$.

Case 2: $|I| < r$. Here, $\{v_i \,:\, i \in I\}$ is linearly independent. If even $\{v_i \,:\, i \in I\} \cup \{w_j \,:\, j \in J\})$ is linearly independent, then the left hand side of \eqref{eq:tangent-stability} equals $\frac{1}{|I|+|J|}(\sum_{i \in I} V_i + |J| W) \le \max(V_0,W)$. Hence, we are necessarily left with the following situation: $\{1, \ldots, s\} \subseteq J$, $0 \in J$, and $a_1 e_1 + \cdots + a_r e_r \in F' := \lin(\{v_i \,:\, i \in I\})$. Now, Lemma~\ref{lemma-criterion} finishes the proof.
\end{proof}

\subsection{Computing the volumes}

Our next goal is make Proposition~\ref{stability-criterion} more explicitly applicable by computing the volumes of these Cayley polytopes.

\begin{proposition}
Let $c_0, \ldots, c_r \in \ZZ_{\ge 0}$, and $\nu > 0$. Then the normalized volume of $(\nu + c_0) \Delta_s * \cdots * (\nu + c_r) \Delta_s$ equals
\[\sum_{k=0}^{s}  \binom{s + r}{k}  \left(\sum_{d_0+\cdots+d_r=s-k} c_0^{d_0} \cdots c_r^{d_r}\right)\nu^k.\]

\label{ausfaktorisiert}
\end{proposition}

In particular, we get:

\begin{corollary}
\[V_0 = \sum_{k=0}^{s}  \binom{s + r-1}{k}  \left(\sum_{d_1+\cdots+d_r=s-k} a_1^{d_1} \cdots a_r^{d_r}\right)\nu^k,\]
%\[W = \sum_{k=0}^{s-1}  \binom{s + r - 1}{k}  \left(\sum_{d_0+\cdots+d_r=s-1-k} a_0^{d_0} \cdots a_r^{d_r}\right)\nu^k.\]
\[W= \sum_{k=0}^{s-1}  \binom{s + r - 1}{k}  \left(\sum_{d_1+\cdots+d_r=s-1-k} a_1^{d_1} \cdots a_r^{d_r}\right)\nu^k,\]
and for $i \in \{1, \ldots, r\}$ 
\[V_i = \sum_{k=0}^{s}  \binom{s + r-1}{k}  \left(\sum_{d_1+\cdits+d_r=s-k} a_1^{d_1} \cdits a_r^{d_r}\right)\nu^k,\]
with the convention that the interior expression equals $0^{s-k}$ if $r=1$ (hence, $V_1 = \nu^s$ in this special case).
\label{vol-formula}
\end{corollary}

In particular, $W \le V_0$, hence, case (c) in Proposition~\ref{stability-criterion} is not necessary to consider.

The proof of Proposition~\ref{ausfaktorisiert} relies on the following general result:

\begin{lemma}
For $(k_0, \ldots, k_r) \in \QQ^{r+1} \setminus \{(0, \ldots, 0)\}$,  the normalized volume of $k_0 \Delta_s * \cdots * k_r \Delta_s$ equals 
\[\sum_{m_0 + \cdots + m_r = s} k_0^{m_0} \cdots k_r^{m_r}\]
where the sum is over $m_0, \ldots, m_r \in \NN$. For $r=1$ this expression gets simplified to $(k_1^{s+1}-k_0^{s+1})/(k_1-k_0)$ if $k_1 \not= k_0$, respectively, to $(s+1) k_0^s$ if $k_1=k_0$.

\label{compute-lemma}
\end{lemma}

\begin{proof} By \cite[6.6]{DaKh} the normalized volume equals $$\sum_{m_0 + \cdots + m_r = s} {\rm MV}((k_0 \Delta_s)^{(m_0)}, \ldots, (k_r \Delta_s)^{(m_r)}),$$ where ${\rm MV}$ denotes the (normalized) mixed volume of $r+1$ many $s$-dimensional polytopes and the exponents indicate the multiplicity of the polytope. Now, the statement follows from multilinearity of the mixed volume and ${\rm MV}(\Delta_s, \ldots, \Delta_s) = 1$. 
\end{proof}

The following useful lemma is straightforward to prove by induction using a well-known identity of binomial coefficients, e.g., \cite[(5.26), table 169]{Concrete}.

\begin{lemma}
Let $d_0, \ldots, d_r, k \in \ZZ_{\ge 0}$, and $r \ge 1$. 
$$\sum_{k_0 + \cdots + k_r = k} \binom{d_0+k_0}{d_0} \cdots \binom{d_r+k_r}{d_r} = \binom{d_0 + \cdots + d_r + r +k}{d_0+\cdots + d_r+r}= \binom{d_0 + \cdots + d_r + r +k}{k}$$
\label{binom-lemma}
\end{lemma}

\begin{proof}[Proof of Proposition~\ref{ausfaktorisiert}]
By Lemma~\ref{compute-lemma}, the normalized volume of $(\nu + c_0) \Delta_s * \cdots * (\nu + c_r) \Delta_s$ equals
\[\sum_{m_0 + \cdots + m_r = s} (\nu+c_0)^{m_0} \cdots (\nu+c_r)^{m_r}\]
\[=\sum_{m_0 + \cdots + m_r = s} \left(\sum_{k_0=0}^{m_0} \binom{m_0}{k_0}\nu^{k_0} c_0^{m_0-k_0}\right) \cdots \left(\sum_{k_r=0}^{m_r} \binom{m_r}{k_r}\nu^{k_r} c_r^{m_r-k_r}\right)\]
\[=\sum_{m_0 + \cdots + m_r = s} \sum_{k=0}^{s} \sum_{k_0+\cdots+k_r=k} \binom{m_0}{k_0} \cdots \binom{m_0}{k_r} \cdot c_0^{m_0-k_0} \cdots c_r^{m_r-k_r} \nu^k\]

\[=\sum_{k=0}^{s} \left(\sum_{d_0+\cdots+d_r=s-k} \left(\sum_{k_0+\cdots+k_r=k} \binom{d_0+k_0}{d_0} \cdots \binom{d_r+k_r}{d_r}\right) \cdot c_0^{d_0} \cdots 
c_r^{d_r}\right)\nu^k.\]
By Lemma~\ref{binom-lemma} this simplifies to
\[\sum_{k=0}^{s} \left(\sum_{d_0+\cdots+d_r=s-k}  \binom{s + r}{k} \cdot c_0^{d_0} \cdots c_r^{d_r}\right)\nu^k.\]
\end{proof}

\subsection{A necessary criterion for stability of \texorpdfstring{$P$}{P}}

Now, we can deduce a strong restriction on the variety.

\begin{proposition}
If $P$ is stable, then $a_r = 1$. If $P$ is semistable, then $a_r=1$ or $s=1$, where in the latter case we have $(a_1, \ldots,a_r) \in \{(0, \ldots, 0, 1), (0, \ldots, 0, 2), (0, \ldots, 0, 1, 1)\}$.\label{bound-prop}
\end{proposition}

\begin{proof}

We abbreviate for $i=1, \ldots, r$, and $k=0, \ldots, s$:
\[D^i_k :=  \sum_{d_1+\cdits+d_r=s-k} a_1^{d_1} \cdits a_r^{d_r},\]
where $D^i_s=1$ (even for $r=1$, see the convention in Corollary~\ref{vol-formula}), and 
\[D^0_k := \sum_{d_1+\cdots+d_r=s-k} a_1^{d_1} \cdots a_r^{d_r},\]
\[W_k := \sum_{d_1+\cdots+d_r=s-1-k} a_1^{d_1} \cdots a_r^{d_r},\]
here, $W_s = 0$.

Let $P$ be stable. By Proposition~\ref{stability-criterion}(a) we get
\[\frac{r+s-1}{r+s} V_0 - \frac{1}{r+s} V_1 - \cdots  - \frac{1}{r+s} V_r - \frac{s+1}{r+s} W< 0.\]
By Corollary~\ref{vol-formula}, this implies
\begin{equation}\sum_{k=0}^{s}  \binom{s + r-1}{k} \alpha_k\, \nu^k < 0,\label{negative}\end{equation}
with
\[\alpha_k = \frac{r+s-1}{r+s} D^0_k - \frac{1}{r+s} D^1_k - \cdots - \frac{1}{r+s} D^r_k - \frac{s+1}{r+s} W_k.\]
Let us assume $a_r \ge 2$. Let $k \in \{0, \ldots, s\}$. As $D^i_k \le D^0_k$ for $i=1, \ldots, r-1$, we get
\[\alpha_k \ge \frac{s}{r+s} D^0_k - \frac{1}{r+s} D^r_k - \frac{s+1}{r+s} W_k\]
We note that 
\[a_r W_k + D^r_k = \left(\sum_{d_1+\cdots+d_r=s-1-k} a_1^{d_1} \cdots a_{r-1}^{d_{r-1}} a_r^{d_r+1}\right) + \left(\sum_{d_1+\cdots+d_{r-1}=s-k} a_1^{d_1} \cdots a_{r-1}^{d_{r-1}}\right) \le D^0_k.\]
Hence, 
\[\alpha_k \ge \frac{s}{r+s} D^0_k - \frac{1}{r+s} D^r_k - \frac{s+1}{a_r (r+s)} (D^0_k - D^r_k) = \frac{(a_r-1)s-1}{a_r(r+s)} D^0_k + \frac{s+1-a_r}{a_r(r+s)}D^r_k.\]
As $a_r \ge 2$, we get $(a_r-1)s-1 \ge 0$, and as $D^0_k \ge D^r_k$, this yields 
\begin{equation}
\alpha_k \ge \frac{(a_r-1)s-1}{a_r(r+s)} D^r_k + \frac{s+1-a_r}{a_r(r+s)}D^r_k = \frac{a_r (s-1)}{a_r(r+s)} D^r_k \ge 0.
\label{negative-2}
\end{equation}
However, this implies
\[\sum_{k=0}^{s}  \binom{s + r-1}{k} \alpha_k\, \nu^k \ge 0,\]
a contradiction to \eqref{negative}.

Now, let $P$ be semistable. In this case, inequality \eqref{negative} becomes
\begin{equation}\sum_{k=0}^{s}  \binom{s + r-1}{k} \alpha_k\, \nu^k \le 0.\label{semi-negative}\end{equation}
Assuming $a_r \ge 2$ and $s \ge 2$, we observe from \eqref{negative-2} that $\alpha_k \ge 0$ for $k=0, \ldots, s-1$, and $\alpha_s > 0$ as $D^r_s=1$. 
Hence, 
\[\sum_{k=0}^{s}  \binom{s + r-1}{k} \alpha_k\, \nu^k \ge \binom{s+r-1}{s} \alpha_s \nu^s >0,\]
a contradiction to \eqref{semi-negative}. Hence, let $s=1$. In this case, inequality \eqref{semi-negative} becomes $\alpha_0 + r \alpha_1 \nu \le 0$ with 
\[\alpha_0 = \frac{r}{r+1}{(a_1 + \cdots + a_r)}-\frac{1}{r+1} (r-1) (a_1 + \cdots + a_r) - \frac{2}{r+1} = \frac{1}{r+1}(a_1 + \cdots + a_r - 2),\]
and $\alpha_1 = \frac{r}{r+1} - \frac{r}{r+1} = 0$. Hence, $a_1 + \cdots + a_r \le 2$ which finishes the proof.

\end{proof}

\subsection{Finishing the proof of Theorem~\ref{thm:kleinschmidt}}
\label{proof-secc}

Let $P$ be stable. By Proposition~\ref{bound-prop} we can restrict ourselves to the following situation:
\[0 = a_1 = \cdots = a_z < a_{z+1} = \cdots = a_r = 1,\]
where $z \in \{0, \ldots, r-1\}$. In this case, Corollary~\ref{vol-formula} yields for $i=0, \ldots, z$
\[V_i = \sum_{k=0}^{s}  \binom{s + r-1}{k}  \left(\sum_{d_{z+1}+\cdots+d_r=s-k} 1 \right)\nu^k\]
\[=\sum_{k=0}^{s}  \binom{s + r-1}{k} \binom{r-z+s-k-1}{s-k} \nu^k;\]
for $i=z+1, \ldots, r$
\[V_i = \sum_{k=0}^{s}  \binom{s + r-1}{k}  \left(\sum_{d_{z+1}+\cdits+d_r=s-k} 1 \cdits 1 \right)\nu^k,\] 
which implies for $z < r-1$
\[V_i= \sum_{k=0}^{s}  \binom{s + r-1}{k}  \binom{r-z-1+s-k-1}{s-k} \nu^k,\]
while we get $V_r = \binom{s + r-1}{s} \nu^s$ if $z=r-1$; finally we have
\[W= \sum_{k=0}^{s-1}  \binom{s + r - 1}{k}  \left(\sum_{d_{z+1}+\cdots+d_r=s-1-k} 1 \right)\nu^k\]
\[= \sum_{k=0}^{s-1}  \binom{s + r - 1}{k}  \binom{r-z+s-1-k-1}{s-1-k}  \nu^k.\]

Let us assume that $P$ is stable and $z < r-1$. Then as above Proposition~\ref{stability-criterion}(a) implies
\[\frac{r+s-1}{r+s} V_0 - \frac{1}{r+s} V_1 - \cdots  - \frac{1}{r+s} V_r - \frac{s+1}{r+s} W< 0.\]
Hence, 
\begin{equation}\sum_{k=0}^{s}  \binom{s + r-1}{k} \alpha_k\, \nu^k < 0,\label{negative2}\end{equation}
where we have for $k=0, \ldots, s-1$
\begin{align*}
	\alpha_k  =&\; \left(\frac{r+s-1}{r+s} - \frac{z}{r+s}\right) \binom{r-z+s-k-1}{s-k} - \frac{r-z}{r+s}  \binom{r-z+s-k-2}{s-k} -
	\\ 
	&\;- \frac{s+1}{r+s} \binom{r-z+s-k-2}{s-1-k} \\
	=&\;\left(\frac{r+s-1-z}{r+s} - \frac{r-z}{r+s}\right)  \binom{r-z+s-k-2}{s-k} + \\
	&\;+ \left(\frac{r+s-1-z}{r+s} - \frac{s+1}{r+s}\right) \binom{r-z+s-k-2}{s-1-k} \\
	=&\;\frac{s-1}{r+s} \binom{r-2-z+s-k}{s-k} + \frac{r-2-z}{r+s} \binom{r-2-z+s-k}{s-1-k} \ge 0,
\end{align*}
and for $k=s$ we have
\begin{equation}\alpha_s=\frac{r+s-1-z}{r+s} -\frac{r-z}{r+s}=\frac{s-1}{r+s} \ge 0,\label{final-s}\end{equation}
a contradiction to \eqref{negative2}.\\

Let $z=r-1$, hence, $W = \sum_{k=0}^{s-1}  \binom{s + r - 1}{k} \nu^k$, $V_r = \binom{s + r-1}{s} \nu^s$, and for $i=0, \ldots, r-1$ we have $V_i =\sum_{k=0}^{s}  \binom{s + r-1}{k} \nu^k$. We have 
\[\frac{1}{r} \sum_{i=0}^r V_i =\left(\sum_{k = 0}^{s-1}   \binom{s+r-1}{k} \nu^{k}\right) + \frac{r+1}{r} \binom{s+r-1}{s} \nu^s > V_0 > W\]
Let us consider the case (d1) in Proposition~\ref{stability-criterion} where we compute for $I \subseteq \{0,\ldots,r\}$ with $r \in I$ and $l := |I| < r$ that $\frac{1}{l+s}(\sum_{i \in I} V_i + (s+1) W)$ equals
\[\frac{1}{l+s}\left((l-1) \left(\sum_{k = 0}^{s}   \binom{s+r-1}{k} \nu^{k}\right) 
+ \binom{s+r-1}{s} \nu^{s}+(s+1) \sum_{k=0}^{s-1} \binom{s+r-1}{k} \nu^k\right)\]
\[=\left(\sum_{k = 0}^{s-1} \binom{s+r-1}{k} \nu^{k}\right) + \frac{l}{l+s} \binom{s+r-1}{s} \nu^{s} < V_0.\]
Finally, let us note that the case (d2) cannot occur as $z=r-1$ and $|I| < r$. Hence, $(\sum_{i=0}^r V_i)/r$ is the maximum of (a)-(d) in Proposition~\ref{stability-criterion}. 

Therefore, Proposition~\ref{stability-criterion} implies that $P$ is stable if and only if 
\[\frac{s}{r(r+s)} \left(\sum_{i=0}^r V_i\right) - \frac{s+1}{r+s} W < 0,\]
equivalently,
\[\frac{s}{r(r+s)} \left(r \left(\sum_{k = 0}^{s}   \binom{s+r-1}{k} \nu^{k}\right) + \binom{s+r-1}{s} \nu^s\right) - \frac{s+1}{r+s} \left(\sum_{k=0}^{s-1}  \binom{s + r - 1}{k} \nu^k\right) < 0,\]
equivalently,
\[p(x)/(r+s) = -\frac{1}{r+s} \left(\sum_{k = 0}^{s-1}  \binom{s+r-1}{k} \nu^{k}\right) + \frac{s(r+1)}{r(r+s)}\binom{s+r-1}{s} \nu^s < 0.\]
We observe that by Descartes' rule of signs there is only one positive root of this polynomial. Note that for $r=1$ we get $p(x) = -\left(\sum_{q = 0}^{s-1} \binom{s}{q} x^{q}\right) + 2s x^{s} = -(x+1)^s + (2s+1) x^s$; 
and for $s=1$ we have $p(x) = -1 + (r+1) x$. 

Finally, let us consider the semistable situation. If $s \ge 2$, then by Proposition~\ref{bound-prop} the same arguments apply to show that we can assume $a_r=1$ and $z=r-1$ (here, the contradiction is that \eqref{negative2} is nonpositive, while \eqref{final-s} is positive). As $(a_1, \ldots,a_r) =(0, \ldots, 0, 1)$, $P$ is then semistable if and only if $p(\nu) \le 0$. So, let $s=1$, where we have to deal with the two remaining cases in Proposition~\ref{bound-prop}. If $(a_1, \ldots,a_r) \in \{(0, \ldots, 0, 2), (0, \ldots, 0, 1, 1)\}$, one computes in both cases for the right hand side in Proposition~\ref{stability-criterion} the value $r \nu+2$, while the value in (b) is equal to $(r+1) \nu + 2$. Hence, $P$ is never semistable in these cases.

%Note that in this Section we have assumed $a_r>0$. If $a_r=0$, then $X$ is a product of projective spaces and the claim is well-known. 
$\hfill \qed$

\section{Two Examples of higher Picard rank}
The following example shows that there exist toric Fano varieties that do not admit a morphism to projective space and whose tangent bundle is stable with respect to the anticanonical polarisation.

\begin{example}\label{ex:NoMapsToPn}
Let $X$ be the $3$-dimensional smooth toric Fano variety that is the blow up at a line of $\mathcal{O}_{\PP^1\times \PP^1}\oplus \mathcal{O}_{\PP^1\times \PP^1}(1,1)$. Then $X$ does not admit a morphism to $\PP^3$ by \cite{WatanabeWatanabe}. However,  the tangent bundle to $X$ is stable with respect to the anticanonical polarisation by \cite{zbMATH00892566}. One can also check this by applying Proposition~\ref{thm:tangent-criterion} to the polytope corresponding to the anticanonical polarisation, whose vertices are \begin{multline*}\{(0,-1,-1), (-1,-1,-1), (-1,0,-1), (0,0,-1), (-1,-1,0), (1,-1,0), (-1,2,1), \\(-1,0,1), (2,0,1), (2,2,1)\}.\end{multline*}
See also \cite[Example 4-11]{Fanography} and 
\cite[ID\#10]{GRDB}. 
\end{example}

The next example shows that the stability region for the tangent bundle inside the nef cone is usually neither convex nor polyhedral. At a first glance this may look like a surprising fact for toric varieties. However, when replacing the self-intersection $(D)^{n-1}$ in the definition of the stability notion by an arbitrary class $\alpha$ of a movable curve, then for a fixed subbundle $\F$ the condition $\mu_\alpha(\F) < \mu_\alpha(\T_X)$ imposes a linear condition on $\alpha$. Since, there are again only finitely many subbundles to consider, these conditions cut out a rational polyhedral subcone of the cone of movable curves (c.f. \cite{Pang}). Then our stability region $\Stab(\T_X) \subset \Nef(X)$ is just the preimage of this polyhedral cone under the non-linear map
\[\Nef(X) \to \Mov(X),\quad [D] \to [(D)^{n-1}].\]
For a systematic treatment of stability with respect to curve classes in the toric setting see \cite{Pang}.

\begin{example}
We consider the iterated blowup $\phi \colon X \to \PP^3$ in a point and the strict transform of a line through this point. Let us denote the pullback of a hyperplane by $H$ and the exceptional divisors of the first and second blowup by $E_1$ and $E_2$ respectively.

The classes of $\CO_X(H)$, $\CO_X(E_1)$ and $\CO_X(E_2)$ form a basis of the Picard group of $X$.  The nef cone is spanned by $\CO_X(H)$, $\CO_{X}(H-E_1)$ and  $\CO_X(H-E_1-E_2)$. Thus, a line bundle of the form 
\[\CO_X(\lambda H +\mu_1(H-E_1) - \mu_2 (E_2))\]
 is ample iff $\lambda, \mu_1, \mu_2 > 0$ and $\mu_1 > \mu_2$. We may rescale such a line bundle and obtain a $\QQ$-line bundle of the form 
\[\CO_X(H +\nu_1(H-E_1) - \nu_2 (E_2)).\]

  To describe the corresponding polytope we consider the standard simplex  $\Delta = \conv \{0,e_1,e_2,e_3\}$ and the following  halfspaces in $\RR^3$. 
  \[H^+_1(\nu)=\{u \in \RR^3 \mid u_1 + u_2 \geq \nu\}, \quad H^+_2(\nu)=\{u \in \RR^3 \mid u_3 \leq \nu\}\]
Then the rational polytope corresponding to our ample line bundle from above is given by
  \[(1+\nu_1)\Delta \cap H^+_1(\nu_1) \cap H^+_2(\nu_2).\]
We have two parallel facets $P^{\rho_0}, P^{-\rho_0} \prec P$, which are perpendicular to $\rho_0 = \RR_{\geq 0}\cdot (0,0,1)$. These facets have volume 
  \(\nu_1^2-\nu_2^2\)
 and
  \((1+\nu_1)^2 - \nu_2^2,\)  respectively.  The remaining facets consist of a rectangle $P^{\rho_1}$ and a trapezoid  $P^{\rho_2}$ opposite to each other and two more trapezoids $P^{\rho_3}$, $P^{\rho_4}$. Here, $\rho_i$ denotes the corresponding rays in the normal fan of $P$.  Elementary calculations shows
  \begin{align*}
    \vol(P^{\rho_1}) &= 2\nu_2\\
    \vol(P^{\rho_2})  &= (1+\nu_1)^2 - \nu_1^2 \\
    \vol(P^{\rho_3})=\vol(P^{\rho_4}) &=  (1+\nu_1)^2 - \nu_1^2 - 2\nu_2.  
  \end{align*}
  
  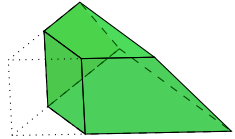
\begin{figure}[ht]
    \centering
    \begin{tikzpicture}
    \rotateRPY{3}{-95}{0}
    \begin{scope}[RPY]
\coordinate (A0) at (0,0,0);
\coordinate (A1) at (-1,0,0);
\coordinate (A2) at (0,0,-1);
\coordinate (A3) at (0,0,-3);
\coordinate (A4) at (-3,0,0);

\coordinate (B0) at (0,1,0);
\coordinate (B1) at (-1,1,0);
\coordinate (B2) at (0,1,-1);
\coordinate (B3) at (0,1,-2);
\coordinate (B4) at (-2,1,0);
\fill [opacity=0.5,fill=green!80!blue] (A4) -- (A3) -- (B3) -- (B4) -- cycle;
\draw [dashed, fill opacity=0.5,fill=green!80!blue] (A4) -- (A1) -- (A2) -- (A3) -- cycle;
\draw [dashed,fill opacity=0.5,fill=green!70!black] (B1) -- (B4) -- (A4) -- (A1) -- cycle;
\draw [fill opacity=0.4,fill=green!80!black] (B4) -- (B1) -- (B2) -- (B3) -- cycle;
\draw [fill opacity=0.4,fill=green!70!black] (B1) --(A1) -- (A2) -- (B2) -- cycle;
\draw [fill opacity=0.4,fill=green!70!black] (B2) -- (B3) -- (A3) -- (A2) -- cycle;
\draw [dotted] (A0) -- (A1) -- (A2) -- cycle;
\draw [dotted] (B0) -- (B1) -- (B2) -- cycle;
\draw [dotted] (A0) -- (B0);
\end{scope}
\end{tikzpicture}
\caption{The polytope $P_D$ corresponding to the polarisation}
\label{fig:F1-bundle}
\end{figure}
We consider the subspaces $F_1=\{\rho_0, -\rho_0\}$ and $F_2=\Span \{\rho_0, -\rho_0, \rho_1, \rho_2\}$. Note, that $\dim F_1=1$ and $\dim F_2 =2$. Now applying the stability codition (\ref{eq:tangent-stability}) for $F_1$ gives.
\begin{equation*}
  \vol P^{\rho_0} + \vol P^{-\rho_0} < \frac{\vol \partial P}{3}
\label{eq:F1}                                          
\end{equation*}
with $\vol \partial P$ being the sum over all facet volumes.
Now, by using the facet volumes stated above, the inequality can be seen to be equivalent to
\begin{equation}
\frac{1}{3} {\left(\nu_{1} + 1\right)}^{2} - \frac{5}{3}  \nu_{1}^{2} + \frac{4}{3} \nu_{2}^{2} - \frac{2}{3} \nu_{2} > 0.\label{eq:F1-2}
\end{equation}
Similarly, evaluating  (\ref{eq:tangent-stability}) for $F_2$ leads to the inequality.
\begin{equation*}
  \vol P^{\rho_0} + \vol P^{-\rho_0} +   \vol P^{\rho_1} +  \vol P^{\rho_2} < \frac{\vol \partial P}{3}
\label{eq:F2}
\end{equation*}
which is equivalent to
\begin{equation}
  \label{eq:F2-2}
\frac{1}{3} {\left(\nu_{1} + 1\right)}^{2} - \frac{2}{3} \nu_{1}^{2} + \frac{1}{3}  \nu_{2}^{2} - \frac{5}{3}  \nu_{2} > 0.
\end{equation}
A sketch of a cross-section of $\Nef(X)$ with the regions cut out by (\ref{eq:F1-2}) and(\ref{eq:F2-2}), respectively, is shown in Figure~\ref{fig:F1-nef}.
\begin{figure}[ht]
  \centering
\resizebox{0.6\textwidth}{!}{\input{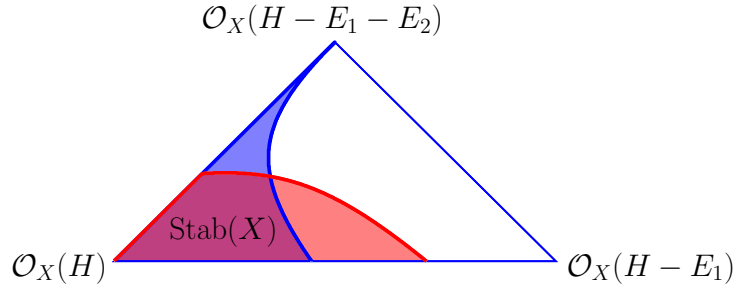}}
\caption{The stability region in $\Nef(X)$}
\label{fig:F1-nef}
\end{figure}
By Remark~\ref{rem:finitely-many-subspaces} there are four more one-dimensional subspaces and
three more two-dimensional subspaces which could provide additional obstructions for the stability of the tangent bundle. However, the corresponding subsheaves turn out to be not destabilising for any ample polarisation.  Hence, $\Stab(\T_X) \subset \Nef(X)$ is cut out by the two inequalities (\ref{eq:F1-2}) and (\ref{eq:F2-2}).
\end{example}
\bibliography{tstab}

\begin{thebibliography}{jBLYZ15}

\bibitem[BDGP18]{BiswasDeyGenvPoddar}
Indranil {Biswas}, Arijit {Dey}, Ozhan {Genc}, and Mainak {Poddar}.
\newblock {On stability of tangent bundle of toric varieties}.
\newblock {\em ArXiv e-prints}, August 2018.

\bibitem[Bel]{Fanography}
Pieter Belmans.
\newblock {Entry in Fanography, A tool to visually study the geography of Fano
  3-folds}.
\newblock
  \href{https://fanography.pythonanywhere.com/4-5}{https://fanography.pythonanywhere.com/4-5}.

\bibitem[BK]{GRDB}
Gavin {Brown} and Alexander {Kasprzyk}.
\newblock {Graded Ring Database. A database of graded rings in algebraic
  geometry}.
\newblock \href{http://www.grdb.co.uk}{http://www.grdb.co.uk}.

\bibitem[CCD97]{CattaniCoxDickenstein97}
Eduardo {Cattani}, David {Cox}, and Alicia {Dickenstein}.
\newblock Residues in toric varieties.
\newblock {\em Compositio Math.}, 108(1):35--76, 1997.

\bibitem[CDS14]{zbMATH06497829}
Xiuxiong {Chen}, Simon {Donaldson}, and Song {Sun}.
\newblock {K\"ahler-Einstein metrics and stability.}
\newblock {\em {Int. Math. Res. Not.}}, 2014(8):2119--2125, 2014.

\bibitem[CDS15a]{zbMATH06394344}
Xiuxiong {Chen}, Simon {Donaldson}, and Song {Sun}.
\newblock {K\"ahler-Einstein metrics on Fano manifolds. I: Approximation of
  metrics with cone singularities.}
\newblock {\em {J. Am. Math. Soc.}}, 28(1):183--197, 2015.

\bibitem[CDS15b]{zbMATH06394345}
Xiuxiong {Chen}, Simon {Donaldson}, and Song {Sun}.
\newblock {K\"ahler-Einstein metrics on Fano manifolds. II: Limits with cone
  angle less than $2\pi$.}
\newblock {\em {J. Am. Math. Soc.}}, 28(1):199--234, 2015.

\bibitem[CDS15c]{zbMATH06394346}
Xiuxiong {Chen}, Simon {Donaldson}, and Song {Sun}.
\newblock {K\"ahler-Einstein metrics on Fano manifolds. III: Limits as cone
  angle approaches $2\pi$ and completion of the main proof.}
\newblock {\em {J. Am. Math. Soc.}}, 28(1):235--278, 2015.

\bibitem[{Dan}78]{zbMATH03661497}
Vladimir~I. {Danilov}.
\newblock {Geometry of toric varieties.}
\newblock {\em {Russ. Math. Surv.}}, 33(2):97--154, 1978.

\bibitem[DDK19]{DDK}
Jyoti {Dasgupta}, Arijit {Dey}, and Bivas {Khan}.
\newblock Stability of equivariant vector bundles over toric varieties.
\newblock Communicated, 2019.

\bibitem[DK86]{DaKh}
Vladimir~I. {Danilov} and A.~G. {Khovanski\u{\i}}.
\newblock Newton polyhedra and an algorithm for calculating {H}odge-{D}eligne
  numbers.
\newblock {\em Izv. Akad. Nauk SSSR Ser. Mat.}, 50(5):925--945, 1986.

\bibitem[{Fah}89]{zbMATH04099463}
Rachid {Fahlaoui}.
\newblock {Stabilit\'e du fibr\'e tangent des surfaces de Del Pezzo. (Stability
  of the tangent bundle of Del Pezzo surfaces).}
\newblock {\em {Math. Ann.}}, 283(1):171--176, 1989.

\bibitem[GKP94]{Concrete}
Ronald~L. {Graham}, Donald~E. {Knuth}, and Oren {Patashnik}.
\newblock {\em Concrete mathematics}.
\newblock Addison-Wesley Publishing Company, Reading, MA, second edition, 1994.
\newblock A foundation for computer science.

\bibitem[GKP16]{zbMATH06561543}
Daniel {Greb}, Stefan {Kebekus}, and Thomas {Peternell}.
\newblock {Movable curves and semistable sheaves.}
\newblock {\em {Int. Math. Res. Not.}}, 2016(2):536--570, 2016.

\bibitem[HL14]{zbMATH06285042}
Martin {Henk} and Eva {Linke}.
\newblock {Cone-volume measures of polytopes.}
\newblock {\em {Adv. Math.}}, 253:50--62, 2014.

\bibitem[jBLYZ13]{zbMATH06168607}
K\'aroly jun. {B\"or\"oczky}, Erwin {Lutwak}, Deane {Yang}, and Gaoyong
  {Zhang}.
\newblock {The logarithmic Minkowski problem.}
\newblock {\em {J. Am. Math. Soc.}}, 26(3):831--852, 2013.

\bibitem[jBLYZ15]{zbMATH06435124}
K\'aroly jun. {B\"or\"oczky}, Erwin {Lutwak}, Deane {Yang}, and Gaoyong
  {Zhang}.
\newblock {Affine images of isotropic measures.}
\newblock {\em {J. Differ. Geom.}}, 99(3):407--442, 2015.

\bibitem[{Kle}88]{zbMATH04085806}
Peter {Kleinschmidt}.
\newblock {A classification of toric varieties with few generators.}
\newblock {\em {Aequationes Math.}}, 35(2-3):254--266, 1988.

\bibitem[Kly90]{Klyachko90}
Alexander~A. Klyachko.
\newblock Equivariant vector bundles on toral varieties.
\newblock {\em Math. USSR-Izv.}, 35(2):337--375, 1990.

\bibitem[{Kly}98]{KlyachkoHermitian}
Alexander~A. {Klyachko}.
\newblock Stable bundles, representation theory and {H}ermitian operators.
\newblock {\em Selecta Math. (N.S.)}, 4(3):419--445, 1998.

\bibitem[{Kob}87]{zbMATH00044936}
Shoshichi {Kobayashi}.
\newblock {\em {Differential geometry of complex vector bundles.}}
\newblock Princeton, NJ: Princeton University Press; Tokyo: Iwanami Shoten
  Publishers, 1987.

\bibitem[{Koo}11]{zbMATH05931673}
Martijn {Kool}.
\newblock {Fixed point loci of moduli spaces of sheaves on toric varieties.}
\newblock {\em {Adv. Math.}}, 227(4):1700--1755, 2011.

\bibitem[KS98]{zbMATH01324687}
Allen {Knutson} and Eric {Sharpe}.
\newblock {Sheaves on toric varieties for physics.}
\newblock {\em {Adv. Theor. Math. Phys.}}, 2(4):873--961, 1998.

\bibitem[LWX18]{li2018algebraicity}
Chi {Li}, Xiaowei {Wang}, and Chenyang {Xu}.
\newblock Algebraicity of the metric tangent cones and equivariant k-stability.
\newblock {\em arXiv:1805.03393}, 2018.

\bibitem[{Lü}83]{zbMATH03890301}
Martin {Lübke}.
\newblock {Stability of Einstein-Hermitian vector bundles.}
\newblock {\em {Manuscr. Math.}}, 42:245--257, 1983.

\bibitem[{Mab}87]{Mabuchi87}
Toshiki {Mabuchi}.
\newblock Einstein-{K}\"{a}hler forms, {F}utaki invariants and convex geometry
  on toric {F}ano varieties.
\newblock {\em Osaka J. Math.}, 24(4):705--737, 1987.

\bibitem[{Mat}57]{zbMATH03149530}
Yoz\^o {Matsushima}.
\newblock {Sur la structure du groupe d'hom\'eomorphismes analytiques d'une
  certaine vari\'et\'e kaehl\'erienne.}
\newblock {\em {Nagoya Math. J.}}, 11:145--150, 1957.

\bibitem[{Oda}88]{zbMATH00192947}
Tadao {Oda}.
\newblock {\em Convex bodies and algebraic geometry}, volume~15 of {\em
  Ergebnisse der Mathematik und ihrer Grenzgebiete (3) [Results in Mathematics
  and Related Areas (3)]}.
\newblock Springer-Verlag, Berlin, 1988.
\newblock An introduction to the theory of toric varieties, Translated from the
  Japanese.

\bibitem[{Pan}15]{Pang}
Thiam-Sun {Pang}.
\newblock {\em The Harder-Narasimhan Filtrations and Rational Contractions}.
\newblock PhD thesis, Universit\"at Freiburg, Fakult\"at f\"ur Mathematik und
  Physik, 2015.

\bibitem[{Per}04]{Perling04}
Markus {Perling}.
\newblock Graded rings and equivariant sheaves on toric varieties.
\newblock {\em Math. Nachr.}, 263/264:181--197, 2004.

\bibitem[{Ste}96]{zbMATH00892566}
Andreas {Steffens}.
\newblock {On the stability of the tangent bundle of Fano manifolds.}
\newblock {\em {Math. Ann.}}, 304(4):635--643, 1996.

\bibitem[WW82]{WatanabeWatanabe}
Keiichi {Watanabe} and Masayuki {Watanabe}.
\newblock The classification of {F}ano {$3$}-folds with torus embeddings.
\newblock {\em Tokyo J. Math.}, 5(1):37--48, 1982.

\bibitem[WZ04]{wang04}
Xu-Jia {Wang} and Xiaohua {Zhu}.
\newblock {K\"ahler--Ricci solitons on toric manifolds with positive first
  Chern class.}
\newblock {\em {Adv. Math.}}, 188(1):87--103, 2004.

\end{thebibliography}
\bibliographystyle{alpha}
\end{document}